\documentclass[12pt]{article}
\usepackage{amssymb}
\usepackage{amsmath}
\usepackage{theorem}

\newtheorem{lem}{Lemma}[section]

\newtheorem{prp}{Proposition}[section]
\theorembodyfont{\rmfamily}

\makeatletter

\@addtoreset{equation}{section}
\makeatother
\def\A{{\text{\boldmath $A$}}}
\def\B{{\text{\boldmath $B$}}}
\def\C{{\text{\boldmath $C$}}}
\def\D{{\text{\boldmath $D$}}}
\def\H{{\text{\boldmath $H$}}}
\def\I{{\text{\boldmath $I$}}}
\def\L{{\text{\boldmath $L$}}}
\def\P{{\text{\boldmath $P$}}}
\def\R{{\text{\boldmath $R$}}}
\def\S{{\text{\boldmath $S$}}}
\def\T{{\text{\boldmath $T$}}}

\def\V{{\text{\boldmath $V$}}}
\def\W{{\text{\boldmath $W$}}}
\def\X{{\text{\boldmath $X$}}}
\def\Y{{\text{\boldmath $Y$}}}
\def\Z{{\text{\boldmath $Z$}}}
\def\la{{\lambda}}
\def\bde{{\text{\boldmath $\delta$}}}
\def\bDe{{\text{\boldmath $\Delta$}}}
\def\bDeh{{\widehat \bDe}}
\def\bGa{{\text{\boldmath $\Gamma$}}}
\def\bOm{{\text{\boldmath $\Omega$}}}
\def\bTh{{\text{\boldmath $\Theta$}}}
\def\bUp{{\text{\boldmath $\Upsilon$}}{}}
\def\bPhi{{\text{\boldmath $\Phi$}}}
\def\bSi{{\text{\boldmath $\Sigma$}}}

\def\Nc{{\cal N}}
\def\Oc{{\cal O}}
\def\Tc{{\cal T}}
\def\Vc{{\cal V}}
\def\Wc{{\cal W}}
\def\Rb{\mathbb{R}}
\def\Db{\mathbb{D}}
\def\zero{{\bf 0}}

\def\diag{{\rm diag}}
\def\tr{{\rm tr}\,}
\def\infi{{\infty}}
\def\non{{\nonumber}}
\newcommand{\qed}{\hfill\hbox{\rule{7pt}{7pt}}}
\def\ia{{\accent 19 \char 16}}
\begin{document}

\title{Estimation of a high-dimensional covariance matrix with the Stein loss}
\author{
Hisayuki Tsukuma\footnote{Faculty of Medicine, Toho University, 
5-21-16 Omori-nishi, Ota-ku, Tokyo 143-8540, Japan, E-Mail: tsukuma@med.toho-u.ac.jp}
}
\maketitle 

\begin{abstract}
The problem of estimating a normal covariance matrix is considered from a decision-theoretic point of view, where the dimension of the covariance matrix is larger than the sample size.
This paper addresses not only the nonsingular case but also the singular case in terms of the covariance matrix.
Based on James and Stein's minimax estimator and on an orthogonally invariant estimator, some classes of estimators are unifiedly defined for any possible ordering on the dimension, the sample size and the rank of the covariance matrix.
Unified dominance results on such classes are provided under a Stein-type entropy loss.
The unified dominance results are applied to improving on an empirical Bayes estimator of a high-dimensional covariance matrix.

\par\vspace{4mm}\noindent
{\it AMS 2010 subject classifications:} Primary 62H12; secondary 62C12.

\par\vspace{2mm}\noindent
{\it Key words and phrases:} Empirical Bayes method, inadmissibility, Moore-Penrose pseudo-inverse, pseudo Wishart distribution, singular multivariate normal distribution, singular Wishart distribution, statistical decision theory.
\end{abstract}

\section{Introduction}
\label{sec:intro}

This paper addresses the problem of a normal covariance matrix relative to the Stein loss, where the dimension of the covariance is larger than the sample size.
This problem is precisely formulated as follows:
Let $\X_1,\X_2,\ldots,\X_n$ be independently and identically distributed as $\Nc_p(\zero_p,\bSi)$.
Assume that $p>n$ and $\bSi$ is a $p\times p$ positive definite matrix of unknown parameters.
Denote $\S=\sum_{i=1}^n\X_i\X_i^t$.
Then $\S$ is distributed as 
\begin{equation}\label{eqn:model}
\S\sim\Wc_p(n,\bSi).
\end{equation}
In the $p>n$ case, Srivastava and Khatri (1979, page 72) and D\ia{}az-Garc\ia{}a et al.\! (1997) called $\Wc_p(n,\bSi)$ the pseudo Wishart distribution with $n$ degrees of freedom and mean $n\bSi$.
We here consider the problem of estimating $\bSi$ relative to the Stein loss
\begin{equation}\label{eqn:loss}
L_p(\bde,\bSi)=\tr\bSi^{-1}\bde -\log\det(\bSi^{-1}\bde)-p,
\end{equation}
where $\bde$ stands for an estimator of $\bSi$.
Assume that, with probability one, $\bde$ is an positive definite matrix based on $\S$.
The accuracy of estimators is measured by the risk function $R_p(\bde,\bSi)=E[L_p(\bde,\bSi)]$, where the expectation is taken with respect to the model (\ref{eqn:model}).

\medskip
If $n\geq p$, then the Wishart matrix $\S$ has the same rank $p$ as the covariance matrix $\bSi$ with probability one.
In such case, many decision-theoretic studies have been done for the problem of estimating $\bSi$ in the literature.
James and Stein (1961) first discussed decision-theoretic estimation of $\bSi$.
They considered the LU decomposition of $\S$ and succeeded to derive a minimax estimator of $\bSi$ relative to the Stein loss (\ref{eqn:loss}).
The James and Stein (1961) minimax estimator, however, depends on the coordinate system.
The dependence results in inadmissibility of their minimax estimator.
Typical improved estimators on James and Stein's minimax estimator are orthogonally invariant estimators, which are not influenced by the coordinate system.
The orthogonally invariant estimators have been proposed by Stein (1975, 1977).
See also Dey and Srinivasan (1985), who gave other dominance results via orthogonally invariant estimators.

\medskip
In the $p>n$ case, Kubokawa and Srivastava (2008) and Konno (2009) studied decision-theoretic covariance estimation relative to quadratic losses.
However, an analytical dominance result in the $p>n$ case with the Stein loss (\ref{eqn:loss}) has not been obtained as yet.

\medskip
This paper gives some dominance results relative to the Stein loss (\ref{eqn:loss}) in the $p>n$ case and extends the dominance results to the case where $\bSi$ is singular.
To this end, Section \ref{sec:unified} starts with unifiedly considering the estimation problem for any possible ordering on $n$, $p$ and the rank of $\bSi$.
The singular case does not allow us to use the Stein loss (\ref{eqn:loss}) because the inverse of the singular $\bSi$ does not exist.
Therefore Section \ref{sec:unified} defines a Stein-like loss function for estimation of the singular $\bSi$.
We give a unified expression of the James and Stein type estimator for all possible orderings on $n$, $p$ and the rank of $\bSi$.
Section \ref{sec:unified} also provides a unified expression of orthogonally invariant estimators improving on the James and Stein type estimator relative to the Stein-like loss.

\medskip
Section \ref{sec:high} mainly discusses the $p>n$ case for estimation of a nonsingular $\bSi$ relative to the usual Stein loss (\ref{eqn:loss}).
An empirical Bayes estimator is derived from an inverse Wishart prior.
Some improving methods on the empirical Bayes estimator are established by using the dominance results obtained in Section \ref{sec:unified}.
The Monte Carlo simulations show that an improved estimator performs well when $p$ is much larger than $n$.
Moreover alternative estimators are unifiedly constructed for both nonsingular and singular cases in terms of $\bSi$.
In Section \ref{sec:remarks}, we give some remarks on our results of this paper and related topics.

\section{Unified dominance results on covariance estimation}
\label{sec:unified}

\subsection{Preliminaries}
\label{subsec:pre}

First, we describe the problem of estimating a covariance matrix unifiedly in the nonsingular and the singular cases.

\medskip
Assume that the $p\times n$ observation matrix $\X$ has the form
\begin{equation}\label{eqn:singular_model}
\X=\B\Z,
\end{equation}
where $\B$ is a $p\times r$ matrix of unknown parameters with $p\geq r$ and $\Z$ is an $r\times n$ random matrix.
Assume that $\B$ is of full column rank, namely $r$, and $r$ is known.
Let all the columns of $\Z$ be independently and identically distributed as $\Nc_r(\zero_r,\I_r)$.
Then the columns of $\X$ are i.i.d.\! sample from $\Nc_p(\zero_p,\bSi)$, where $\bSi=\B\B^t$ is a positive semi-definite matrix of rank $r$.
Denote
$$
\S=\X\X^t,
$$
which follows $\Wc_p(n,\bSi)$.
In the case where $r<p$, $\Nc_p(\zero_p,\bSi)$ and $\Wc_p(n,\bSi)$ represent, respectively, the singular multivariate normal and the singular Wishart distributions.
For the definition of the singular distributions, see Srivastava and Khatri (1979, pages 43 and 72) and also D\ia{}az-Garc\ia{}a et al.\! (1997).
Note also that $\bSi$ is of rank $r$, while $\S$ is of rank $\min(n,r)$ with probability one.

\medskip
In this section, we handle only estimators which are positive semi-definite matrices of rank
$$
q=\min(n,r)
$$
with probability one.
Write such estimators as $\bde_q$.
Moreover, $\bde_q$ are also assumed to satisfy the condition that the rank of $\bSi^+\bde_q$ is $q$ with probability one, where $\bSi^+$ is the Moore-Penrose pseudo-inverse of $\bSi$.
Since $\bde_q$ and $\bSi^+$ are positive semi-definite, the $q$ nonzero eigenvalues of $\bSi^+\bde_q$ are positive.
Note that $\tr\bSi^+\bde_q$ is equal to a sum of all the positive eigenvalues of $\bSi^+\bde_q$.
Both nonsingular and singular cases of the Stein loss (\ref{eqn:loss}) are unifiedly defined as
\begin{equation}\label{eqn:singular_loss}
L_q(\bde_q,\bSi)=\tr\bSi^+\bde_q-\log\pi(\bSi^+\bde_q)-q,
\end{equation}
where $\pi(\bSi^+\bde_q)$ stands for a product of all the positive eigenvalues of $\bSi^+\bde_q$.
Then we consider the problem of estimating $\bSi$ relative to the Stein loss (\ref{eqn:singular_loss}).
The corresponding risk function is denoted by
\begin{equation}\label{eqn:singular_risk}
R_q(\bde_q,\bSi)=E[L_q(\bde_q,\bSi)],
\end{equation}
where the expectation is taken with respect to the model (\ref{eqn:singular_model}).

\medskip
Next, we define some notation.
Let $\Oc(r)$ be the group of orthogonal matrices of order $r$.
For $p\geq r$, the Stiefel manifold is denoted by $\Vc_{p,r}=\{\A\in\Rb^{p\times r}:\A^t\A=\I_r\}$.
It is noted that $\Vc_{r,r}=\Oc(r)$.
Let $\Db_r$ be a set of $r\times r$ diagonal matrices whose diagonal elements $d_1,\ldots,d_r$ satisfy $d_1>\cdots>d_r>0$.
Denote by $\Tc_q^+$ the group of lower triangular matrices with positive diagonal elements.

\medskip
The Stein loss (\ref{eqn:singular_loss}) depends on the Moore-Penrose pseudo-inverse of \bSi.
Here some properties are listed for the Moore-Penrose pseudo-inverse.
The proof of the following lemma is given in Harville (1997, Chapter 20).
\begin{lem}\label{lem:mpp}
Let $\B$ be a $p\times r$ matrix of full column rank.
Then the Moore-Penrose pseudo-inverse $\B^+$ of $\B$ uniquely exists and has the following properties:
\newcounter{mpp}
\begin{list}{}{
\topsep=0pt
\parsep=0pt
\parskip=4pt
\itemsep=4pt
\itemindent=0mm
\labelwidth=8pt
\labelsep=6pt
\leftmargin=26pt
\listparindent=12pt
\usecounter{mpp}
}
\renewcommand{\makelabel}{\rm(\arabic{mpp})}

\item $\B^+=(\B^t\B)^{-1}\B^t$;
\item $\H^+=\H^t$ for $\H\in\Vc_{p,r}$;
\item $\B^+=\B^{-1}$ for a nonsingular matrix $\B$;
\item $(\B^+)^t=(\B^t)^+$;
\item $(\B\C^t)^+=(\C^t)^+\B^+$ for a $q\times r$ matrix $\C$ of full column rank.
\end{list}
\end{lem}

\subsection{Constant multiple estimators}
\label{subsec:cme}

Consider a simple class of estimators whose forms are a constant multiple of $\S$.
The simple class is represented by
\begin{equation}\label{eqn:cme}
\bde_q^C(a)=a\S,
\end{equation}
where $a$ is a positive constant and $q=\min(n,r)$.
This class includes the unbiased estimator of $\bSi$,
$$
\bde_q^{UB}=\frac{1}{n}\S.
$$

However $\bde_q^{UB}$ is not the best estimator among the class (\ref{eqn:cme}) relative to the Stein loss (\ref{eqn:singular_loss}).
Note by Lemma \ref{lem:mpp} that $\bSi^+=(\B\B^t)^+=(\B^t)^+\B^+$ and
$$
\bSi^+\S=(\B^t)^+\B^+\B\Z\Z^t\B^t=(\B^t)^+\Z\Z^t\B^t,
$$
which implies that $\bSi^+\bde_q^C(a)$ has the same rank as $\Z\Z^t$.
\begin{prp}\label{prp:bc}
Define $m=\max(n,r)$ and
$$
a_m=\frac{1}{m}.
$$
Then $\bde_q^{BC}=\bde_q^C(a_m)$ is the best estimator among the class {\rm(\ref{eqn:cme})} relative to the Stein loss {\rm(\ref{eqn:singular_loss})}.
Hence for $r>n$, $\bde_q^{BC}$ dominates $\bde_q^{UB}$ relative to the Stein loss {\rm(\ref{eqn:singular_loss})}.
\end{prp}

{\bf Proof.}\ \ 
The nonzero eigenvalues of $\bSi^+\S$ are identical to those of $\B^+\S(\B^t)^+$, so that the nonzero eigenvalues of $\bSi^+\S$ are identical to those of the full-rank matrix
$$
\begin{cases}
\Z\Z^t & \textup{for $n\geq r$},\\
\Z^t\Z & \textup{for $n< r$}.\\
\end{cases}
$$
Since the number of nonzero eigenvalues of $\bSi^+\S$ is $q=\min(n,r)$ with probability one, we obtain $\pi(a\bSi^+\S)=a^q\pi(\Z\Z^t)$, so that the risk of $\bde_q^C(a)$ with respect to the Stein loss (\ref{eqn:singular_loss}) is expressed as
\begin{align*}
R_q(\bde_q^C(a),\bSi)
&=na\,\tr\bSi^+\bSi -q\log a -E[\log\pi(\Z\Z^t)]-q \\
&=nra-q\log a -E[\log\pi(\Z\Z^t)]-q.
\end{align*}
The risk of $\bde_q^C(a)$ is minimized by $\bde_q^C(a_m)$ with
$$
a_m=\frac{q}{nr}=\frac{1}{m}.
$$
Thus the proof is complete.
\qed

\medskip
It follows from equation (82) of James and Stein (1961) that
\begin{equation}\label{eqn:digamma}
E[\log\pi(\Z\Z^t)]=\sum_{i=1}^qE[\log s_i],
\end{equation}
where $s_i\sim\chi^2_{m-i+1}$.
Hence $\bde_q^{BC}$ has the constant risk
\begin{equation}\label{eqn:risk_bc}
R_q(\bde_q^{BC},\bSi)=q\log m - \sum_{i=1}^qE[\log s_i].
\end{equation}

\subsection{The James and Stein type estimator}
\label{subsec:js}

We next construct a James and Stein (1961) like estimator of $\bSi$ for any possible ordering on $n$, $p$ and $r$.

\medskip
Using the same arguments as in Srivastava (2003, equation (2.2)), we can write the $p\times n$ random matrix $\X$ as a block matrix
$$
\begin{pmatrix}
\X_{11} & \X_{12} \\
\X_{21} & \X_{22}
\end{pmatrix}
,
$$
where $\X_{11}$ is a $q\times q$ nonsingular matrix.
Recall that $\X=\B\Z$. 
Partition $\B$ and $\Z$ into block matrices as, respectively,
$$
\B=\begin{pmatrix} \B_1 \\ \B_2 \end{pmatrix},\qquad
\Z=(\Z_1,\Z_2),
$$
where $\B_1$ and $\Z_1$ are, respectively, $q\times r$ and $r\times q$ matrices.
Note that $\X_{11}=\B_1\Z_1$.
Since $\X_{11}$ is nonsingular, $\B_1$ has a full row rank.
Thus, there exist unique elements $\bTh\in\Tc_q^+$ and $\bGa_1\in\Vc_{r,q}$ such that $\B_1=\bTh\bGa_1^t$.
The decomposition $\B_1^t=\bGa_1\bTh^t$ represents the QR decomposition of $\B_1^t$.
Also, $\B_1=\bTh\bGa_1^t$ is called the LQ decomposition of $\B_1$.
For the uniqueness of the QR decomposition, see Harville (1997, page 67).

\medskip
Take $\bGa_2\in\Vc_{r,r-q}$ such that $\bGa=(\bGa_1,\bGa_2)\in\Oc(r)$.
For $r\leq n$, the LQ decomposition of $\bGa^t\Z$ can be written as $\Y\V^t$, where $\Y\in\Tc_r^+$ and $\V\in\Vc_{n,r}$.
For $r>n$, $\bGa^t\Z$ is denoted by the block matrix
$$
\bGa^t\Z=\begin{pmatrix} \Z_1 \\ \Z_2 \end{pmatrix},
$$
where $\Z_1$ and $\Z_2$ are, respectively, $n\times n$ and $(r-n)\times n$ matrices.
The LQ decomposition of $\Z_1$ can be written as $\Y_1\V^t$ for $\Y_1\in\Tc_n^+$ and $\V\in\Oc(n)$, which gives that
$$
\bGa^t\Z=\begin{pmatrix} \Y_1 \\ \Z_2\V \end{pmatrix}\V^t=\begin{pmatrix} \Y_1 \\ \Y_2 \end{pmatrix}\V^t,
$$
where $\Y_2=\Z_2\V$.
Hence $\bGa^t\Z$ can uniquely and unifiedly be expressed as
\begin{equation}\label{eqn:GZ}
\bGa^t\Z=\Y\V^t\qquad \textup{for}\ \ \V\in\Vc_{n,q}\ \ \textup{and}\ \ 
\Y=\begin{pmatrix} \Y_1 \\ \Y_2 \end{pmatrix},
\end{equation}
where $\Y_1\in\Tc_q^+$ and $\Y_2$ is a $(r-q)\times q$ matrix.

\medskip
Let $\C=\B\bGa$, which is written as
\begin{equation}\label{eqn:C}
\C=\begin{pmatrix} \bTh\bGa_1^t \\ \B_2 \end{pmatrix}(\bGa_1,\bGa_2)
=\begin{pmatrix}
\bTh & \zero_{q\times (r-q)} \\
\B_2\bGa_1 & \B_2\bGa_2 \end{pmatrix}
.
\end{equation}
Combining (\ref{eqn:GZ}) and (\ref{eqn:C}), we can uniquely decompose $\X$ as
\begin{equation}\label{eqn:LQ}
\X=\B\bGa\bGa^t\Z=\T\V^t,
\end{equation}
where
$$
\T=\C\Y
=\begin{pmatrix}
\bTh\Y_1 \\
\B_2\bGa_1\Y_1+ \B_2\bGa_2\Y_2
\end{pmatrix}
.
$$
It is then noted that $\bTh\Y_1\in\Tc_q^+$.

\medskip
The probability distributions of nonzero elements of $\Y$ are given as follows.
\begin{lem}\label{lem:Ldis}
For $i=1,\ldots,q$ and $j=i,\ldots,r$, denote by $y_{j,i}$ the $(j,i)$-th element of $\Y$.
Then all the elements $y_{j,i}$'s are mutually independent and
$$
y_{i,i}^2\sim \chi^2_{n-i+1},\qquad y_{j,i}\sim\Nc(0,1)\qquad (i=1,\ldots, q,\ j=i+1,\ldots,r).
$$
\end{lem}

{\bf Proof.}\ \ 
It is noted that $\bGa^t\Z\sim\Nc_{r\times n}(\zero_{r\times n},\I_r\otimes\I_n)$.
For the $n\geq r$ and $n<r$ cases, see Lemma 3.2.1 of Srivastava and Khatri (1979) and Corollary 3.1 of Srivastava (2003), respectively.
\hfill$\Box$

\medskip
Applying (\ref{eqn:LQ}) to the Wishart matrix $\S=\X\X^t$ gives that
$$
\X\X^t=\T\T^t=\begin{pmatrix} \T_1 \\ \T_2 \end{pmatrix} ( \T_1^t, \T_2^t),
$$
where $\T=( \T_1^t, \T_2^t)^t$ is $p\times q$ matrix such that $\T_1=\bTh\Y_1\in\Tc_q^+$ and $\T_2$ is a $(p-q)\times q$ matrix.
Then we consider the class of estimators, which has the form
\begin{equation}\label{eqn:class_tr}
\bde_q^T=\T\D_q\T^t,
\end{equation}
where $\D_q=\diag(d_1,\ldots,d_q)$ and the $d_i$'s are positive constants.

\begin{prp}\label{prp:js}
Let $\D_q^{JS}=\diag(d_1^{JS},\ldots,d_q^{JS})$, where $d_i^{JS}=(n+r-2i+1)^{-1}$ for $i=1,\ldots,q$.
Then the best estimator among the class {\rm(\ref{eqn:class_tr})} relative to the loss {\rm(\ref{eqn:singular_loss})} is given by $\bde_q^{JS}=\T\D_q^{JS}\T^t$, which dominates $\bde_q^{BC}$.
\end{prp}

{\bf Proof.}\ \ 
Noting from Lemma \ref{lem:mpp} that
$$
\C^t\bSi^+\C=\bGa^t\B^t(\B\B^t)^+\B\bGa=\bGa^t\B^t(\B^t)^+\B^+\B\bGa=\I_r,
$$
we observe 
\begin{equation}\label{eqn:tr_Si_deT}
E[\tr\bSi^+\bde_q^T]=E[\tr\D_q\T^t\bSi^+\T]=E[\tr\D_q\Y^t\C^t\bSi^+\C\Y]=E[\tr\D_q\Y^t\Y].
\end{equation}

\medskip
In the $p\geq r>n$ case, we partition $\Y$ as $\Y=(\Y_1^t,\Y_2^t)^t$, where $\Y_1\in\Tc_n^+$.
From Lemma \ref{lem:Ldis}, it follows that $E[\Y_2^t\Y_2]=(r-n)\I_n$ and $E[\Y_1^t\Y_1]$ is the diagonal matrix of order $n$ with the $i$-th diagonal element
$$
\sum_{j=i}^nE[y_{j,i}^2]=E[y_{i,i}^2]+\sum_{j>i}^nE[y_{j,i}^2]=(n-i+1)+(n-i)=2n-2i+1.
$$
Thus we obtain
\begin{equation}\label{eqn:r1_n<r}
E[\tr\D_n\Y^t\Y]=E[\tr\D_n\Y_1^t\Y_1+\tr\D_n\Y_2^t\Y_2]=\sum_{i=1}^n (n+r-2i+1)d_i.
\end{equation}
When $n\geq r$, it follows that
\begin{equation}\label{eqn:r1_n>r}
E[\tr\D_r\Y^t\Y]=\sum_{i=1}^r \sum_{j=i}^r E[d_i y_{j,i}^2]=\sum_{i=1}^r (n+r-2i+1)d_i.
\end{equation}
Combining (\ref{eqn:tr_Si_deT}), (\ref{eqn:r1_n<r}) and (\ref{eqn:r1_n>r}) gives that
\begin{equation}\label{eqn:r1}
E[\tr\bSi^+\bde_q^T]=\sum_{i=1}^q (n+r-2i+1)d_i.
\end{equation}
It is seen that $\bSi^+\bde_q^T$ has the same nonzero eigenvalues as $\D_q\T^t\bSi^+\T$, which implies that
$$
\pi(\bSi^+\bde_q^T)
=\pi(\D_q\T^t\bSi^+\T)\\
=\pi(\D_q\Y^t\Y)
$$
Since $\Y^t\Y$ is a $q\times q$ square matrix of full rank, it follows that
\begin{align}
\pi(\bSi^+\bde_q^T)&=\det(\D_q\Y^t\Y) \non\\
&=\det(\D_q)\det(\Y^t\Y) \non\\
&=\det(\D_q)\pi(\Z\Z^t).
\label{eqn:r2}
\end{align}
Using (\ref{eqn:r1}) and (\ref{eqn:r2}), we can write the risk of $\bde_q^T$ under the loss (\ref{eqn:singular_loss}) as
$$
R_q(\bde_q^T,\bSi)=\sum_{i=1}^q \{(n+r-2i+1)d_i-\log d_i\}-E[\log\pi(\Z\Z^t)]-q.
$$
Hence the $d_i$'s minimizing the risk are given by $d_i^{JS}=(n+r-2i+1)^{-1}$ for $i=1,\ldots,q$.

\medskip
Since $\bde_q^{BC}$ belongs to the class (\ref{eqn:class_tr}), $\bde_q^{JS}$ dominates $\bde_q^{BC}$ relative to the Stein loss (\ref{eqn:singular_loss}).
In fact, $\bde_q^{JS}$ has the constant risk
\begin{equation}\label{eqn:risk_JS}
R_q(\bde_q^{JS},\bSi)=\sum_{i=1}^q \log(n+r-2i+1)-E[\log\pi(\Z\Z^t)],
\end{equation}
which implies by (\ref{eqn:digamma}) and (\ref{eqn:risk_bc}) that
$$
R_q(\bde_q^{JS},\bSi)-R_q(\bde_q^{BC},\bSi)=\sum_{i=1}^q \log(n+r-2i+1)-q\log m<0,
$$
where the inequality follows from concavity of the logarithmic function.
Thus the proof is complete.
\qed

\medskip
The probability density function of $\T$ can be derived explicitly.
The $n\geq p=r$ case is obtained from, for example, Srivastava and Khatri (1979, Lemma 3.2.2). 
For the $p>r$ case, see Srivastava (2003, Theorem 5.2) and D\ia{}az-Garc\ia{}a and Gonz\'alez-Far\ia{}as (2005, Corollary 4).

\subsection{Orthogonally invariant estimators}
\label{subsec:orth}

Make the QR decomposition of $\B$ into $\bUp\B_0^t$, where $\bUp\in\Vc_{p,r}$ and $\B_0\in\Tc_r^+$.
We can uniquely express $\S$ as $\S=\bUp\B_0^t\Z\Z^t\B_0\bUp^t$.
Define $\W=\B_0^t\Z\Z^t\B_0$, which is distributed as $\Wc_r(n,\bOm)$ with $\bOm=\B_0^t\B_0$, where $\bOm$ is positive definite.
The eigenvalue decomposition of $\W$ is written as $\R\L\R^t$, where $\L\in\Db_q$, $\R\in\Vc_{r,q}$ and $q=\min(n,r)$.
Hence we can decompose $\S$ as
\begin{equation*}
\S=\H\L\H^t, 
\end{equation*}
where $\L\in\Db_q$ and $\H=\bUp\R\in\Vc_{p,q}$.

\medskip
Consider the class of estimators
$$
\bde_q^O=\bde_q^O(\S)=\H\bPhi(\L)\H^t,
$$
where $\bPhi(\L)=\diag(\phi_1(\L),\ldots,\phi_q(\L))$ and the $\phi_i(\L)$'s are absolutely continuous functions of $\L$.
The class $\bde_q^O$ is orthogonally invariant in the sense that it satisfies $\O\bde_q^O(\S)\O^t=\bde_q^O(\O\S\O^t)$ for any $\O\in\Oc(p)$.

\medskip
To evaluate the risk of $\bde_q^O$, we require the following lemma.
\begin{lem}\label{lem:st_id}
Abbreviate $\phi_i(\L)$ by $\phi_i$.
Denote $\L=\diag(\ell_1,\ldots,\ell_q)$.
Then we have
$$
E[\tr\bSi^+\H\bPhi(\L)\H^t]=E\bigg[\sum_{i=1}^q\bigg\{(|n-r|-1)\frac{\phi_i}{\ell_i}+2\frac{\partial\phi_i}{\partial\ell_i}+2\sum_{j>i}^q\frac{\phi_i-\phi_j}{\ell_i-\ell_j}\bigg\}\bigg].
$$
\end{lem}

{\bf Proof.}\ \ 
It follows from Lemma \ref{lem:mpp} that
$$
\bSi^+=(\bUp\bOm\bUp^t)^+=(\bOm\bUp^t)^+\bUp^+=(\bUp^t)^+\bOm^+\bUp^+=\bUp\bOm^{-1}\bUp^t, 
$$
so that
\begin{equation}\label{eqn:trace}
E[\tr\bSi^+\H\bPhi(\L)\H^t]=E[\tr\bOm^{-1}\R\bPhi(\L)\R^t].
\end{equation}
Recall that $\W\sim\Wc_r(n,\bOm)$ and the eigenvalue decomposition of $\W$ is denoted by $\R\L\R^t$.
The remainder of the proof for $n\geq r$ is based on the same arguments as in Sheena (1995, Section 2) and, for $n<r$, on those as in Kubokawa and Srivastava (2008, Lemma A.1).
Their results are applied to the r.h.s.\! of (\ref{eqn:trace}), so we get this lemma.
\qed

\begin{lem}\label{lem:risk_orth}
The risk of $\bde_q^O$ under the loss {\rm (\ref{eqn:singular_loss})} is written as
\begin{align*}
R(\bde_q^O,\bSi)
&=E\bigg[\sum_{i=1}^q\bigg\{(|n-r|-1)\frac{\phi_i}{\ell_i}+2\frac{\partial\phi_i}{\partial\ell_i}+2\sum_{j>i}^q\frac{\phi_i-\phi_j}{\ell_i-\ell_j}-\log\frac{\phi_i}{\ell_i}\bigg\}\bigg] \\
&\qquad -E[\log\pi(\Z\Z^t)]-q.
\end{align*}
\end{lem}

{\bf Proof.}\ \ 
Note that 
\begin{equation}\label{eqn:pi_orth}
\pi(\bSi^+\bde_q^O)=\pi(\Z\Z^t)\det(\L^{-1}\bPhi(\L)).
\end{equation}
Using (\ref{eqn:pi_orth}) and Lemma \ref{lem:st_id} gives the risk expression of this lemma.
\qed

\medskip
The following proposition results from Lemma \ref{lem:risk_orth}.
\begin{prp}
Define
$$
\bde_q^{ST}=\H\L\D_q^{JS}\H^t.
$$
Then $\bde_q^{ST}$ dominates $\bde_q^{JS}$ relative to the Stein loss {\rm (\ref{eqn:singular_loss})}.
\end{prp}

{\bf Proof.}\ \ 
We can prove this proposition in the same way as in Dey and Srinivasan (1985, Theorem 3.1).
Using Lemma \ref{lem:risk_orth} and (\ref{eqn:risk_JS}), we can write the difference in risk of $\bde_q^{ST}$ and $\bde_q^{JS}$ as $R_q(\bde_q^{ST},\bSi)-R_q(\bde_q^{JS},\bSi)=E[\bDeh^{ST}]$, where
$$
\bDeh^{ST}=\sum_{i=1}^q\bigg\{(|n-r|+1)d_i^{JS}+2\sum_{j>i}^q\frac{d_i^{JS}\ell_i-d_j^{JS}\ell_j}{\ell_i-\ell_j}\bigg\}-q.
$$
Hence if $\bDeh^{ST}\leq 0$, then $\bde_q^{ST}$ dominates $\bde_q^{JS}$.
It is observed that
\begin{align*}
\sum_{i=1}^q\sum_{j>i}^q\frac{d_i^{JS}\ell_i-d_j^{JS}\ell_j}{\ell_i-\ell_j}
&=\sum_{i=1}^q\sum_{j>i}^q\frac{d_i^{JS}(\ell_i-\ell_j)+(d_i^{JS}-d_j^{JS})\ell_j}{\ell_i-\ell_j} \\
&<\sum_{i=1}^q\sum_{j>i}^q d_i^{JS}
=\sum_{i=1}^q(q-i) d_i^{JS},
\end{align*}
where the inequality is verified by the ordering properties $\ell_1>\cdots>\ell_q$ and $d_1^{JS}<\cdots<d_q^{JS}$.
Thus we obtain
\begin{align*}
\bDeh^{ST}
<\sum_{i=1}^q(|n-r|+1+2q-2i)d_i^{JS}-q 
=\sum_{i=1}^q(n+r-2i+1)d_i^{JS}-q=0,
\end{align*}
which completes the proof.
\qed

\medskip
Besides $\bde_q^{ST}$ given above, many types of orthogonally invariant estimators are proposed for the $n\geq p=r$ case.
See, for example, Stein (1977), Dey and Srinivasan (1985), Haff (1991), Perron (1992), Sheena and Takemura (1992) and Yang and Berger (1994).
Their results would be applicable to the cases when $\min(n,p)\geq r$ and $p\geq r>n$.

\section{Estimation of a high-dimensional covariance matrix}
\label{sec:high}

\subsection{An empirical Bayes estimator}
\label{subsec:EB}

We here deal with the problem of estimating $\bSi$ in the model (\ref{eqn:model}) relative to the usual Stein loss (\ref{eqn:loss}).
Note that the covariance matrix $\bSi$ is of rank $p$, while the Wishart matrix $\S$ is of rank $n$.
Using an empirical Bayes method, we first provide a full-rank estimator as a target which should be improved.

\medskip
Denote $\X=(\X_1,\ldots,\X_n)$, where the $\X_i$'s are i.i.d.\! sample from $\Nc_p(\zero_p,\bSi)$.
Note that $\X$ is a $p\times n$ matrix and $\S=\X\X^t$.
Then the likelihood of $\bSi$ is proportional to
$$
L(\bSi|\X)\propto (\det\bSi)^{-n/2}\exp\Big(-\frac{1}{2}\tr\bSi^{-1}\X\X^t\Big).
$$
Assume that $\bSi$ has a prior density 
\begin{equation*}
p(\bSi|\la)\propto (\det\bSi)^{-(k+p+1)/2}\exp\Big(-\frac{\la}{2}\tr\bSi^{-1}\Big),\quad \la>0.
\end{equation*}
The resulting Bayes estimator $\bde_p^{Bayes}$ is written as
\begin{equation}\label{eqn:Bayes}
\bde_p^{Bayes}=\frac{1}{n+k}(\X\X^t+\la\I_p)=\frac{1}{n+k}(\S+\la\I_p).
\end{equation}
Here we estimate $\la$ from the marginal density of $\X$,
\begin{align*}
p(\X|\la)
&=K \la^{kp/2}\{\det(\X\X^t+\la\I_p)\}^{-(n+k)/2}\\
&=K \la^{kp/2}\{\det(\X^t\X+\la\I_n)\}^{-(n+k)/2},
\end{align*}
where $K$ is a normalizing constant.
Since $\det(\X^t\X+\la\I_n)=\prod_{j=1}^n(\ell_j+\la)$, where the $\ell_j$'s are eigenvalues of $\X^t\X$, the logarithm of the marginal density $p(\X|\la)$ is given by
\begin{equation*}
\log p(\X|\la)= \frac{kp}{2}\log\la-\frac{n+k}{2}\sum_{j=1}^n \log(\ell_j+\la)+\log K,
\end{equation*}
which is used to obtain
\begin{equation*}
\frac{\partial}{\partial\la}\log p(\X|\la) = \frac{kp}{2}\la^{-1}-\frac{n+k}{2}\sum_{j=1}^n \frac{1}{\ell_j+\la}=0,
\end{equation*}
namely, the maximum likelihood estimator of $\la$ is a solution of
\begin{equation*}
\sum_{j=1}^n \frac{\la}{\ell_j+\la}=\frac{kp}{n+k}.
\end{equation*}
Denote by $\hat{\la}^{ML}$ the resulting maximum likelihood estimator of $\la$.
Substitute $\hat{\la}^{ML}$ for $\la$ in (\ref{eqn:Bayes}), we get the empirical Bayes estimator
\begin{equation}\label{eqn:B}
\bde_p^B=\frac{1}{n+k} (\S+\hat{\la}^{ML}\I_p ).
\end{equation}

\medskip
Motivated by (\ref{eqn:B}) and taking account of Proposition \ref{prp:bc}, we define the class of estimators as
\begin{equation}\label{eqn:EB}
\bde_p^{EB}(b)=a_p(\S+\hat{\la}_b\I_p),
\end{equation}
where $a_p=p^{-1}$, $b=b(\S)$ is a differentiable bounded function of $\S$, and $\hat{\la}_b\ (\geq 0)$ satisfies
\begin{equation}\label{eqn:la}
\sum_{j=1}^n\frac{\hat{\la}_b}{\ell_j+\hat{\la}_b}=b.
\end{equation}
For existence of a unique solution $\hat{\la}_b$, $b$ requires at least that $0\leq b<n$.
Note also that $\bde_p^{EB}(b)$ is of full-rank with probability one.

\medskip
To compare risk functions, we need the lower and the upper bounds of $\hat{\la}_b$.
Note from Lemma \ref{lem:mpp} that 
\begin{equation*}
\sum_{i=1}^n\ell_j^{-1}=\tr(\X^t\X)^{-1}=\tr\X(\X^t\X)^{-2}\X^t=\tr(\X\X^t)^+=\tr\S^+.
\end{equation*}
Also, note that $\sum_{j=1}^n\ell_j=\tr\S$.
\begin{lem}\label{lem:bound_la}
The lower and the upper bounds of $\hat{\la}_b$ are given as follows.
$$
\frac{b}{n-b}\cdot\frac{n}{\tr\S^+}\leq \hat{\la}_b\leq\frac{b}{n-b}\cdot\frac{\tr\S}{n}.
$$
When $b<1$, it particularly holds that $\hat{\la}_b\leq b/\{(1-b)\tr\S^+\}$.
\end{lem}

{\bf Proof.}\ \ 
Let $f(x|c)=c/(x+c)$ for a positive constant $c$.
Since $f(x|c)$ is convex in $x$ for $x\geq 0$, it is observed that
$$
b=\sum_{j=1}^nf(\ell_j|\hat{\la}_b)\geq n f\bigg(\frac{1}{n}\sum_{j=1}^n \ell_j\bigg|\hat{\la}_b\bigg)
=\frac{n^2\hat{\la}_b}{\tr\S+n\hat{\la}_b},
$$
which gives the upper bound of $\hat{\la}_b$.
Next, let $g(x)=x/(1+x)$ for $x\geq 0$.
The concavity of $g(x)$ leads to
$$
b=\sum_{j=1}^n\frac{\hat{\la}_b\ell_j^{-1}}{1+\hat{\la}_b\ell_j^{-1}}=\sum_{j=1}^ng(\hat{\la}_b\ell_j^{-1})\leq n g\bigg(\frac{1}{n}\sum_{j=1}^n \hat{\la}_b\ell_j^{-1}\bigg)
=\frac{n\hat{\la}_b\tr\S^+}{n+\hat{\la}_b\tr\S^+},
$$
which gives the lower bound of $\hat{\la}_b$.

\medskip
When $b<1$, we can see that
\begin{align*}
b&=\sum_{j=1}^n\frac{\hat{\la}_b\ell_j^{-1}}{1+\hat{\la}_b\ell_j^{-1}}
\geq \sum_{j=1}^n\frac{\hat{\la}_b\ell_j^{-1}}{1+\hat{\la}_b\sum_{j=1}^n\ell_j^{-1}}=\frac{\hat{\la}_b\tr\S^+}{1+\hat{\la}_b\tr\S^+},
\end{align*}
which yields that $\hat{\la}_b\leq b/\{(1-b)\tr\S^+\}$.
\qed

\medskip
The finiteness of the risk of $\bde_p^{EB}(b)$ is verified in the following lemma.
\begin{lem}\label{lem:finite_HF}
Assume that there exist positive constants $B_1$ and $B_2$ such that $B_1\leq b\leq B_2<n$.
If $p-n-1>0$, then the risk of $\bde_p^{EB}(b)$ is finite.
\end{lem}

{\bf Proof.}\ \ 
A simple calculation yields that
\begin{align*}
R_p(\bde_p^{EB}(b),\bSi)
&=E\left[a_p\hat{\la}_b\tr\bSi^{-1}-\sum_{i=1}^n\log(\ell_i+\hat{\la}_b)-(p-n)\log\hat{\la}_b\right]\\
&\qquad\qquad +npa_p-p\log a_p+\log\det\bSi-p .
\end{align*}
From the given assumption, there exist positive constants $C_1$ and $C_2$ such that $C_1/n\leq b/(n-b)\leq n C_2$.
Using Lemma \ref{lem:bound_la}, we observe that
$$
\log C_1-\log\tr\S^+\leq \log\hat{\la}_b\leq \log C_2+\log\tr\S.
$$
The well-known inequalities $1-x^{-1}\leq\log x\leq x-1$ for $x>0$ imply that $E[\log\hat{\la}_b]$ is finite if $E[\tr\S^+]<\infi$.
Under the same condition, we can verify the finiteness of $E[\hat{\la}_b]$ and $E[\sum_{i=1}^n\log(\ell_i+\hat{\la}_b)]$.

\medskip
Note that $E[\tr\S^+]=E[\tr(\X^t\X)^{-1}]=E[\tr(\Z^t\bSi\Z)^{-1}]$, where $\Z\sim\Nc_{p\times n}(\zero_{p\times n},\I_p\otimes\I_n)$, so that
$$
0<E[\tr\S^+]<E[\tr(\Z^t\Z)^{-1}]\tr\bSi^{-1}.
$$
Note also that $\Z^t\Z\sim\Wc_n(p,\I_n)$.
Thus for $p-n-1>0$, it follows that $E[\tr(\Z^t\Z)^{-1}]=n(p-n-1)^{-1}$, which completes the proof.
\qed

\subsection{Dominance results}
\label{subsec:dominance}

In the case that $p=r>n$, define the eigenvalue decomposition of $\S$ as $\S=\H\L\H^t$ with $\H\in\Vc_{p,n}$ and $\L=\diag(\ell_1,\ldots,\ell_n)\in\Db_n$.
Take $\H_0$ as a $p\times (p-n)$ matrix such that $(\H,\H_0)\in\Oc(p)$.
Consider here the following shrinkage estimator
\begin{align*}
\bde_p^{SH}(b)
&=\bde_p^{EB}(b)-a_p\hat{\la}_b\H\H^t \\
&=a_p(\S+\hat{\la}_b\H_0\H_0^t),
\end{align*}
where $\hat{\la}_b$ and $b$ are defined in (\ref{eqn:EB}).
The rank of $\bde_p^{SH}(b)$ is $p$ with probability one.
If $p-n-1>0$, the risk of $\bde_p^{SH}(b)$ is finite, which is verified in the same way as Lemma \ref{lem:finite_HF}.
The following proposition can be obtained for domination of $\bde_p^{SH}(b)$ over $\bde_p^{EB}(b)$.
\begin{prp}\label{prp:SH}
In the model {\rm(\ref{eqn:model})}, we consider the problem of estimating $\bSi$ relative to the usual Stein loss {\rm(\ref{eqn:loss})}.
Assume that there exists a positive constant $C$ such that $b\leq C<n$.
Let $c_0=6(n+1)/(3p-4n-4)$ for $3p-4n-4>0$.
If $b\geq c_0n/(1+c_0) $ and $\sum_{i=1}^n\partial b/\partial \ell_i\geq 0$, then $\bde_p^{SH}(b)$ dominates $\bde_p^{EB}(b)$ relative to the usual Stein loss {\rm(\ref{eqn:loss})}.
\end{prp}

The proof of Proposition \ref{prp:SH} requires suitable bounds of the logarithmic function $\log(1+x)$.
Here we employ an upper and a lower bounds of $\log(1+x)$ based on the Pad\'e approximants.
For details of the Pad\'e approximants, see Baker and Graves-Morris (1996).
The approximants yield the following simple lemma, whose proof is omitted since it can easily be verified.
\begin{lem}\label{lem:pade}
For $x\geq 0$, it follows that
$$
\frac{2x}{2+x}\leq \log(1+x) \leq \frac{x(6+x)}{2(3+2x)}.
$$
The upper and the lower bounds given above are concave in $x$.
\end{lem}

\medskip
{\bf Proof of Proposition \ref{prp:SH}.}\ \ 
Note that
$$
\bde_p^{SH}(b)=a_p\{\H\L\H^t+\H_0(\hat{\la}_b\I_{p-n})\H_0^t\}
$$
and also 
$$
\bde_p^{EB}(b)
=a_p\{\H(\L+\hat{\la}_b\I_n)\H^t+\H_0(\hat{\la}_b\I_{p-n})\H_0^t\}.
$$
The difference in risk of $\bde_p^{SH}(b)$ and $\bde_p^{EB}(b)$ is written by
\begin{align}\label{eqn:rd_m-eb}
&R_p(\bde_p^{SH}(b),\bSi)-R_p(\bde_p^{EB}(b),\bSi) \non \\
&=E[-a_p\hat{\la}_b\tr\bSi^{-1}\H\H^t-\log\det\L+\log\det(\L+\hat{\la}_b\I_n)] \non \\
&=E\bigg[-a_p\hat{\la}_b\tr\bSi^{-1}\H\H^t+\sum_{i=1}^n \log(1+\hat{\la}_b\ell_i^{-1})\bigg].
\end{align}
Using Lemma \ref{lem:st_id} gives that $R_p(\bde_p^{SH}(b),\bSi)-R_p(\bde_p^{EB}(b),\bSi)=E[\bDeh{}^{SH}]$, where
\begin{equation}\label{eqn:m-eb1}
\bDeh{}^{SH}=\sum_{i=1}^n\bigg\{-a_p(p-n-1)\hat{\la}_b\ell_i^{-1}-2a_p\frac{\partial \hat{\la}_b}{\partial \ell_i}+ \log(1+\hat{\la}_b\ell_i^{-1})\bigg\}.
\end{equation}
Thus, if $\bDeh{}^{SH}\leq 0$, then $\bde_p^{SH}(b)$ dominates $\bde_p^{EB}(b)$.

\medskip
Differentiating both sides of (\ref{eqn:la}) with respect to $\ell_i$ yields that 
$$
\bigg(\frac{\partial \hat{\la}_b}{\partial \ell_i}\bigg)\sum_{j=1}^n\frac{1}{\ell_j+\hat{\la}_b}
-\frac{\hat{\la}_b}{(\ell_i+\hat{\la}_b)^2}-\bigg(\frac{\partial \hat{\la}_b}{\partial \ell_i}\bigg)\sum_{j=1}^n\frac{\hat{\la}_b}{(\ell_j+\hat{\la}_b)^2}=\frac{\partial b}{\partial \ell_i},
$$
so that
\begin{equation}\label{eqn:diff_la}
\sum_{i=1}^n\frac{\partial \hat{\la}_b}{\partial \ell_i}=\frac{\sum_{i=1}^n\hat{\la}_b(\ell_i+\hat{\la}_b)^{-2}+\sum_{i=1}^n\partial b/\partial\ell_i}{\sum_{i=1}^n\ell_i(\ell_i+\hat{\la}_b)^{-2}}\geq 0.
\end{equation}

\medskip
Let $f(x)=x(6+x)/(6+4x)$.
Using Lemma \ref{lem:pade}, we observe that
\begin{align}\label{eqn:m-eb2}
\sum_{i=1}^n \log(1+\hat{\la}_b\ell_i^{-1})
&\leq \sum_{i=1}^n f(\hat{\la}_b\ell_i^{-1}) \non \\
&\leq nf\bigg(\frac{1}{n}\sum_{i=1}^n \hat{\la}_b\ell_i^{-1}\bigg)
=\frac{\hat{\la}_b\tr\S^+(6n+\hat{\la}_b\tr\S^+)}{6n+4\hat{\la}_b\tr\S^+},
\end{align}
where the second inequality follows from concavity of $f(x)$.
Combining (\ref{eqn:rd_m-eb}), (\ref{eqn:m-eb1}) and (\ref{eqn:m-eb2}) gives that
\begin{align*}
\bDeh{}^{SH}
&\leq -a_p(p-n-1)\hat{\la}_b\tr\S^+ +\frac{\hat{\la}_b\tr\S^+(6n+\hat{\la}_b\tr\S^+)}{6n+4\hat{\la}_b\tr\S^+} \\
&=a_p\hat{\la}_b\tr\S^+\times \frac{6n(n+1)-(3p-4n-4)\hat{\la}_b\tr\S^+}{6n+4\hat{\la}_b\tr\S^+}.
\end{align*}
Using the lower bound of $\hat{\la}_b$ given in Lemma \ref{lem:bound_la}, we can see that $\bDeh{}^{SH}\leq 0$ if $3p-4n-4>0$ and
$$
6n(n+1)-(3p-4n-4)\frac{bn}{n-b}\leq 0,
$$
namely $b\geq c_0n/(1+c_0)$.
Hence the proof is complete.
\qed

\medskip
We give two examples for $b$.
First, $b$ is restricted to a positive constant.
The estimator $\bde_p^{SH}(b)$ can be written as
\begin{equation*}
\bde_p^{SH}(b)=\bde_n^{BC}+a_p\hat{\la}_b\H_0\H_0^t,
\end{equation*}
where $\bde_n^{BC}$ is given by (\ref{eqn:cme}).
The risk of $\bde_p^{SH}(b)$ can alternatively be expressed as
\begin{equation}\label{eqn:risk_M}
R_p(\bde_p^{SH}(b),\bSi)=R_n(\bde_n^{BC},\bSi)+R_{p-n}(a_p\hat{\la}_b\H_0\H_0^t,\bSi),
\end{equation}
where $R_n$ and $R_{p-n}$ are defined in (\ref{eqn:singular_risk}).
It is much hard to find out an optimal constant for $b$.
Furthermore, the performance of $\bde_p^{SH}(b)$ would worsen if $b$ is too large.
So we take
\begin{equation}\label{eqn:b_0}
b_0=\frac{c_0n}{1+c_0}.
\end{equation}
The resulting estimator $\bde_p^{SH}(b_0)$ dominates $\bde_p^{EB}(b_0)$ when $3p-4n-4>0$.

\medskip
Next, consider
\begin{equation}\label{eqn:b_1}
b_1=b_1(\S)=(1+\ell_n/\ell_1)b_0.
\end{equation}
Note that $\ell_1\geq \ell_n$, so $b_1$ is bounded below and above as $b_0\leq b_1\leq 2b_0$.
Also, it is observed that
$$
\sum_{i=1}^n\frac{\partial b_1}{\partial\ell_i}
=b_0\ell_1^{-2}(\ell_1-\ell_n) \geq 0.
$$
Hence it is seen from Proposition \ref{prp:SH} that $\bde_p^{SH}(b_1)$ dominates $\bde_p^{EB}(b_1)$ relative to the usual Stein loss (\ref{eqn:loss}) for $3p-4n-4>0$.

\medskip
The risk expression (\ref{eqn:risk_M}) suggests further modified estimators
\begin{equation*}
\bde_p^{mJS}(b)=\bde_n^{JS}+a_p\hat{\la}_b\H_0\H_0^t,
\end{equation*}
and
\begin{equation}\label{eqn:mST}
\bde_p^{mST}(b)=\bde_n^{ST}+a_p\hat{\la}_b\H_0\H_0^t,
\end{equation}
where $\bde_n^{JS}$ and $\bde_n^{ST}$ are defined in Subsections \ref{subsec:js} and \ref{subsec:orth}, respectively.
Then the following proposition can be proved in the same way as in Subsections \ref{subsec:js} and \ref{subsec:orth}.
\begin{prp}
In the model {\rm(\ref{eqn:model})}, we consider the problem of estimating $\bSi$ relative to the usual Stein loss {\rm(\ref{eqn:loss})}.
Under the assumptions of Proposition \ref{prp:SH}, $\bde_p^{SH}(b)$ is dominated by $\bde_p^{mJS}(b)$, and moreover $\bde_p^{mJS}(b)$ is dominated by $\bde_p^{mST}(b)$.
\end{prp}

\medskip
Proposition \ref{prp:SH} suggests that $\bde_p^{SH}(b)$ dominates $\bde_p^{EB}(b)$ if they depend on a common large $b$.
For a small $b$, it seems to hold the reverse dominance relation.
In fact, we obtain the following proposition.
\begin{prp}\label{prp:SH*}
In the model {\rm(\ref{eqn:model})}, we consider the problem of estimating $\bSi$ relative to the usual Stein loss {\rm(\ref{eqn:loss})}.
Assume that $n\geq 2$ and $p-n-1>0$.
Let $c_*=2(n-1)/(p-n+1)$ and 
\begin{equation}\label{eqn:b_*}
b_*=c_*/(1+c_*).
\end{equation}
If $C_0\leq b\leq b_*$ for a positive constant $C_0$ and $\sum_{i=1}^n\partial b/\partial \ell_i\leq 0$, then $\bde_p^{EB}(b)$ dominates $\bde_p^{SH}(b)$ relative to the usual Stein loss {\rm(\ref{eqn:loss})}.
\end{prp}

{\bf Proof.}\ \ 
In the similar way to (\ref{eqn:diff_la}), it is seen that
\begin{align}\label{eqn:diff_la*}
\sum_{i=1}^n\frac{\partial \hat{\la}_b}{\partial \ell_i}
&\leq \frac{\sum_{i=1}^n\hat{\la}_b(\ell_i+\hat{\la}_b)^{-2}}{\sum_{i=1}^n\ell_i(\ell_i+\hat{\la}_b)^{-2}}
=\sum_{i=1}^n\frac{\hat{\la}_b}{\ell_i}\times \frac{\ell_i(\ell_i+\hat{\la}_b)^{-2}}{\sum_{j=1}^n\ell_j(\ell_j+\hat{\la}_b)^{-2}} \non \\
&\leq \hat{\la}_b\tr\S^+.
\end{align}
It follows from Lemma \ref{lem:pade} that
\begin{align}\label{eqn:m-eb2*}
\sum_{i=1}^n \log(1+\hat{\la}_b\ell_i^{-1})
&\geq \sum_{i=1}^n \frac{2\hat{\la}_b\ell_i^{-1}}{2+\hat{\la}_b\ell_i^{-1}} \non \\
&\geq \sum_{i=1}^n \frac{2\hat{\la}_b\ell_i^{-1}}{2+\hat{\la}_b\sum_{j=1}^n\ell_j^{-1}}
=\frac{2\hat{\la}_b\tr\S^+}{2+\hat{\la}_b\tr\S^+}.
\end{align}
Combining (\ref{eqn:m-eb1}), (\ref{eqn:diff_la*}) and (\ref{eqn:m-eb2*}) gives that
\begin{align*}
\bDeh{}^{SH}
&\geq -a_p(p-n-1)\hat{\la}_b\tr\S^+ -2a_p\hat{\la}_b\tr\S^++\frac{2\hat{\la}_b\tr\S^+}{2+\hat{\la}_b\tr\S^+} \\
&=a_p\hat{\la}_b\tr\S^+\times \frac{2(n-1)-(p-n+1)\hat{\la}_b\tr\S^+}{2+\hat{\la}_b\tr\S^+}.
\end{align*}
If $b\leq b_*\ (<1)$, using the upper bound of Lemma \ref{lem:bound_la} for $b<1$ leads to
\begin{align*}
2(n-1)-(p-n+1)\hat{\la}_b\tr\S^+
&\geq 2(n-1)-(p-n+1)\frac{b}{1-b} \\
&\geq 2(n-1)-(p-n+1)\frac{b_*}{1-b_*}=0,
\end{align*}
which implies that $\bDeh{}^{SH}\geq 0$.
Thus the proof is complete.
\qed

\medskip
Assume that $b$ is a small constant satisfying $b\leq b_*$.
The estimator $\bde_p^{EB}(b)$ is expressed as
$$
\bde_p^{EB}(b)=a_p(\S+\hat{\la}_{b}\I_p)
=a_p\{\H(\L+\hat{\la}_{b}\I_n)\H^t+\H_0(\hat{\la}_{b}\I_{p-n})\H_0^t\},
$$
so the $(p-n)$ eigenvalues among the $p$ nonzero eigenvalues of $\bde_p^{EB}(b)$ are identically $a_p\hat{\la}_{b}=\hat{\la}_{b}/p$.
It is seen from Lemma \ref{lem:bound_la} that
$$
\frac{b}{n-b}\cdot\frac{n}{\tr\S^+}\leq \hat{\la}_{b}\leq \frac{b}{1-b}\cdot\frac{1}{\tr\S^+}.
$$
Note that $n\ell_n^{-1}\geq\tr\S^+\geq\ell_n^{-1}$, so that $\ell_n/n\leq(\tr\S^+)^{-1}\leq\ell_n$.
Moreover in the large-$p$ and small-$n$ case, $c_*$ and $b_*$ probably is a very small value.
Then $\hat{\la}_{b}/p$ may become extremely small, which implies that $\bde_p^{EB}(b)$ may loss stability and deteriorate in performance.
Therefore from Proposition \ref{prp:SH*}, $b_*$ may be a better choice for $b$.
See the next subsection, which gives some simulated values of the risk of $\bde_p^{EB}(b_*)$.

\medskip
We can treat the Haff (1980) type empirical Bayes estimator
$$
\bde_p^{HF}(c)=a_p(\S+cu\I_p),
$$
where $u=1/\tr\S^+$ and $c$ is a positive constant.
Some dominance results on $\bde_p^{HF}(c)$ and $\bde_p^{*SH}(c)=\bde_p^{HF}-cu\H\H^t$ can be derived, and the details are omitted.

\subsection{Monte Carlo studies}\label{subsec:monte}

The Monte Carlo experiments have been performed for comparing the risks of some estimators for some $p$ and $n$.
Each experiment is based on 2,000 independent replications.
We have investigated estimators $\bde_p^{EB}(b)$ and $\bde_p^{mST}(b)$, which are defined in (\ref{eqn:EB}) and (\ref{eqn:mST}), respectively.
It has been assumed that $b=b_0$ and $b_1$, which are given in (\ref{eqn:b_0}) and (\ref{eqn:b_1}), respectively.
Also the risk of $\bde_p^{EB}(b_*)$ has been estimated in our experiments, where $b_*$ is given by (\ref{eqn:b_*}).

\medskip
Note that $b_0$, $b_1$ and $b_*$ satisfy $b(\S)=b(c\S)$ for any positive number $c$.
Also, when $\S$ is transformed into $c\S$ for a positive number $c$, $\hat{\la}_{b}$ satisfying $b(\S)=b(c\S)$ becomes $c\hat{\la}_{b}$.
Hence the risks of $\bde_p^{EB}(b)$ and $\bde_p^{mST}(b)$ with $b=b_0$, $b_1$ and $b_*$ are invariant under the scale transformation $\S\to c\S$ and $\bSi\to c\bSi$ for any positive number $c$.
Furthermore the risks of $\bde_p^{EB}(b)$ and $\bde_p^{mST}(b)$ are invariant under the orthogonal transformation $\S\to\P\S\P^t$ and $\bSi\to \P\bSi\P^t$ for any $\P\in\Oc(p)$.

\medskip
In our experiments, it has been assumed, without loss of generality, that $\bSi$ is a diagonal matrix whose diagonal elements (namely, eigenvalues) are larger than or equal to one.
The following diagonal matrices were considered for an unknown covariance $\bSi$ which should be estimated:
\begin{enumerate}
\setlength{\itemsep}{4pt}
\setlength{\parskip}{0pt}

\item[1)] $\I_p$;
\item[2)] $\diag(10,10^{1-1/p},10^{1-2/p},\ldots,10^{1-(p-2)/p},10^{1-(p-1)/p})$;
\item[3)] $\diag(100,100^{1-1/p},100^{1-2/p},\ldots,100^{1-(p-2)/p},100^{1-(p-1)/p})$.
\end{enumerate}
In Case 1), all the eigenvalues of $\bSi$ are identical.
In Case 2) and $3)$, the eigenvalues of $\bSi$ are widely scattered and the largest eigenvalue is about tenfold or hundredfold of the smallest eigenvalue.

\begin{table}[htb]
\caption{Simulated risk with respect to the usual Stein loss.}
\label{tab:1}
\vspace{3pt}
\centering{{\small
$
\begin{array}{crrc@{\hspace{12pt}}r@{\hspace{8pt}}r@{\hspace{16pt}}r@{\hspace{8pt}}r@{\hspace{20pt}}r}
\hline
\bSi&\multicolumn{1}{c}{p}&\multicolumn{1}{c}{n}&
    &\multicolumn{1}{c}{\bde_p^{EB}(b_0)}
    &\multicolumn{1}{c}{\bde_p^{mST}(b_0)}\quad
    &\multicolumn{1}{c}{\bde_p^{EB}(b_1)}
    &\multicolumn{1}{c}{\bde_p^{mST}(b_1)}\quad
    &\multicolumn{1}{c}{\bde_p^{EB}(b_*)} \\
\hline
1)& 50& 5&& 28.6\,(0.07)& 28.5\,(0.08)& 18.4\,(0.08)& 18.2\,(0.09)&113.4\,(0.10)\\
  &   &15&&  5.0\,(0.01)&  2.7\,(0.02)&  4.8\,(0.01)&  2.0\,(0.02)& 86.2\,(0.06)\\
  &   &25&& 24.3\,(0.08)&  8.7\,(0.02)& 26.5\,(0.10)&  9.6\,(0.03)& 67.3\,(0.06)\\[1pt]
  &100& 5&&115.8\,(0.13)&115.8\,(0.13)& 82.4\,(0.16)& 82.3\,(0.16)&301.1\,(0.14)\\
  &   &25&& 13.3\,(0.02)& 10.8\,(0.03)& 10.4\,(0.02)&  7.3\,(0.03)&228.7\,(0.07)\\
  &   &50&& 41.2\,(0.08)& 14.1\,(0.02)& 44.3\,(0.09)& 15.4\,(0.02)&165.3\,(0.06)\\[1pt]
  &150& 5&&230.7\,(0.16)&230.7\,(0.16)&170.9\,(0.22)&170.9\,(0.22)&516.7\,(0.18)\\
  &   &40&& 18.0\,(0.02)& 13.7\,(0.03)& 14.9\,(0.01)&  9.6\,(0.03)&377.5\,(0.06)\\
  &   &75&& 59.0\,(0.07)& 19.9\,(0.02)& 63.1\,(0.08)& 21.5\,(0.02)&276.1\,(0.06)\\[1pt]
2)& 50& 5&& 26.0\,(0.07)& 25.9\,(0.08)& 19.1\,(0.07)& 18.9\,(0.08)&105.6\,(0.12)\\
  &   &15&& 13.3\,(0.02)& 11.0\,(0.01)& 14.4\,(0.03)& 11.6\,(0.02)& 81.7\,(0.07)\\
  &   &25&& 44.2\,(0.13)& 27.6\,(0.06)& 46.3\,(0.14)& 28.9\,(0.07)& 65.1\,(0.06)\\[1pt]
  &100& 5&&103.0\,(0.14)&102.9\,(0.14)& 75.4\,(0.17)& 75.3\,(0.17)&282.9\,(0.17)\\
  &   &25&& 23.7\,(0.01)& 21.3\,(0.01)& 23.7\,(0.01)& 20.8\,(0.01)&217.6\,(0.08)\\
  &   &50&& 76.7\,(0.12)& 48.2\,(0.06)& 79.5\,(0.13)& 49.9\,(0.06)&160.5\,(0.06)\\[1pt]
  &150& 5&&207.4\,(0.18)&207.4\,(0.18)&155.4\,(0.23)&155.4\,(0.23)&488.0\,(0.21)\\
  &   &40&& 35.6\,(0.01)& 31.3\,(0.01)& 36.1\,(0.01)& 31.1\,(0.01)&361.2\,(0.08)\\
  &   &75&&110.7\,(0.11)& 69.6\,(0.06)&114.5\,(0.12)& 71.9\,(0.06)&268.6\,(0.06)\\[1pt]
3)& 50& 5&& 35.7\,(0.02)& 35.6\,(0.02)& 38.4\,(0.08)& 38.2\,(0.07)& 87.1\,(0.14)\\
  &   &15&& 61.9\,(0.17)& 59.2\,(0.16)& 65.1\,(0.20)& 62.2\,(0.18)& 70.8\,(0.09)\\
  &   &25&&125.6\,(0.33)&105.7\,(0.27)&127.5\,(0.34)&107.2\,(0.28)& 59.8\,(0.07)\\[1pt]
  &100& 5&& 87.4\,(0.11)& 87.4\,(0.11)& 77.2\,(0.08)& 77.2\,(0.08)&237.3\,(0.21)\\
  &   &25&& 99.4\,(0.12)& 96.7\,(0.12)&104.7\,(0.15)&101.8\,(0.14)&188.8\,(0.10)\\
  &   &50&&219.5\,(0.31)&186.2\,(0.25)&221.8\,(0.32)&188.1\,(0.26)&148.2\,(0.07)\\[1pt]
  &150& 5&&165.1\,(0.18)&165.1\,(0.18)&135.9\,(0.17)&135.9\,(0.17)&415.0\,(0.26)\\
  &   &40&&154.4\,(0.13)&149.7\,(0.12)&161.3\,(0.16)&156.1\,(0.14)&318.2\,(0.10)\\
  &   &75&&317.2\,(0.31)&269.7\,(0.25)&320.3\,(0.31)&272.2\,(0.26)&249.2\,(0.07)\\
\hline
\end{array}
$
}}
\end{table}

\medskip
Table \ref{tab:1} shows some simulated risk values.
In the table, the value in parentheses stands for estimated standard error on risk.
For reference, the exact risk of James and Stein's (1961) minimax estimator are $37.096$ ($p=n=50$), $72.0995$ ($p=n=100$) and $106.959$ ($p=n=150$), which can be computed from (\ref{eqn:risk_JS}) and (\ref{eqn:digamma}) of this paper.

\medskip
For large $n\ (=p/2)$, $\bde_p^{mST}(b)$ provides substantial reduction in risk of $\bde_p^{EB}(b)$, but almost not for small $n\ (=5)$.
In the large-$n$ case, $\bde_p^{mST}(b_0)$ is slightly better than $\bde_p^{mST}(b_1)$ and, in the small-$n$ case, $\bde_p^{mST}(b_1)$ is the best estimator among estimators considered here.

\medskip
The estimator $\bde_p^{EB}(b_*)$ has an undesirable performance when $n=5$, and however it enhances the performance as $n$ increases for each $p$.
In Case 3) with large $n\ (=p/2)$, $\bde_p^{EB}(b_*)$ has the smallest risk.

\medskip
All the risks of estimators investigated here considerably varies with the change of $p$, $n$ and $\bSi$.
For example, the risks of $\bde_p^{EB}(b_0)$ and $\bde_p^{EB}(b_*)$ have very different behavior with increasing $n$.
Our numerical results suggest that, although an optimal selection of $b$ would involve difficulty in practical application, we could recommend $\bde_p^{mST}(b_1)$ if $p$ is much larger than $n$.

\subsection{A unified dominance result including both nonsingular and singular cases}

In Subsection \ref{subsec:dominance}, we provided some dominance results for $p=r>n$.
The dominance results can be extended to all possible cases of orderings on $n$, $p$ and $r$ in the model (\ref{eqn:singular_model}).

\medskip
Note that the possible orderings on $n$, $p$ and $r$ are expressed by either $\min(n, p)\geq r$ or $p\geq r>n$.
Let $q=\min(n,r)$ and $m=\max(n,r)$.
The eigenvalue decomposition of $\S$ is written as $\H\L\H^t$, where $\H\in\Vc_{p,q}$ and $\L=\diag(\ell_1,\ldots,\ell_q)\in\Db_q$.
Take $\H_0\in\Vc_{p,p-q}$ such that $(\H,\H_0)\in\Oc(p)$.
Let $\hat{\la}_b$ be a unique solution of the equation
$$
\sum_{i=1}^q\frac{\la}{\ell_i+\la}=b,
$$
where $b$ is a differentiable function of $\S$ and satisfies $0\leq b<q$.
The estimators $\bde_r^{EB}(b)$ and $\bde_r^{SH}(b)$ are defined by, respectively,
\begin{equation}\label{eqn:EB_all}
\begin{split}
\bde_r^{EB}(b)=\begin{cases}
a_n(\S+\hat{\la}_b\H\H^t) & \textup{for $\min(n, p)\geq r$}, \\
a_r(\S+\hat{\la}_b\I_p), & \textup{for $p\geq r>n$},
\end{cases}
\end{split}
\end{equation}
\begin{equation}\label{eqn:SH_all}
\begin{split}
\bde_r^{SH}(b)=\bde_r^{EB}(b)-a_m\hat{\la}_b\H\H^t=\begin{cases}
a_n\S & \textup{for $\min(n, p)\geq r$}, \\
a_r(\S+\hat{\la}_b\H_0\H_0^t), & \textup{for $p\geq r>n$},
\end{cases}
\end{split}
\end{equation}
where $a_m=m^{-1}$.

\medskip
In the $\min(n,p)\geq r$ and the $p=r>n$ cases, the definition (\ref{eqn:EB_all}) and (\ref{eqn:SH_all}) imply that $\bde_r^{EB}(b)$ and $\bde_r^{SH}(b)$ have the same rank as $\bSi$.
However, in the $p>r>n$ case, $\bde_r^{EB}(b)$ and $\bde_r^{SH}(b)$ are of rank $p$, while $\bSi^+\bde_r^{EB}(b)$ and $\bSi^+\bde_r^{SH}(b)$ are of rank $r$.
This is verified as follows:
When $p>r>n$, recall that $\H=\bUp\R$ where $\bUp\in\Vc_{p,r}$ and $\R\in\Vc_{r,n}$, which are defined in the beginning of Subsection \ref{subsec:orth}.
Take $\bUp_0\in\Vc_{p,p-r}$ and $\R_0\in\Vc_{r,r-n}$ such that $(\bUp,\bUp_0)\in\Oc(p)$ and $(\R,\R_0)\in\Oc(r)$.
Define $\H_0\H_0^t=\bUp\R_0\R_0^t\bUp^t+\bUp_0\bUp_0^t$.
Then it is seen that
$$
\H\H^t+\H_0\H_0^t=\bUp(\R\R^t+\R_0\R_0^t)\bUp^t+\bUp_0\bUp_0^t
=\bUp\bUp^t+\bUp_0\bUp_0^t=\I_p
$$
and 
$$
\bUp^t\H_0\H_0^t\bUp=\R_0\R_0^t.
$$
Since $\bSi^+=\bUp\bOm^{-1}\bUp^t$, where $\bUp\in\Vc_{p,r}$ and $\bOm$ is $r\times r$ positive definite, it is observed that
\begin{align*}
\bUp^t\bde_r^{EB}(b)\bUp&=a_r(\R(\L+\hat{\la}_b\I_n)\R^t+\hat{\la}_b\R_0\R_0^t), \\
\bUp^t\bde_r^{SH}(b)\bUp&=a_r(\R\L\R^t+\hat{\la}_b\R_0\R_0^t),
\end{align*}
so that $\bde_r^{EB}(b)$ and $\bde_r^{SH}(b)$ are of rank $p$, while $\bSi^+\bde_r^{EB}(b)$ and $\bSi^+\bde_r^{SH}(b)$ are of rank $r$.
In such $p>r>n$ case, $\bUp\R_0\R_0^t\bUp^t$ is not observable.
Thus it is hard to find an estimator $\bde$ satisfying that both $\bde$ and $\bSi^+\bde$ are of rank $r$.

\medskip
The difference in risk of $\bde_r^{SH}(b)$ and $\bde_r^{EB}(b)$ with respect to the Stein loss (\ref{eqn:singular_loss}) can be written as
\begin{align*}
R_r(\bde_r^{SH}(b),\bSi)-R_r(\bde_r^{EB}(b),\bSi)=
E[-a_m\hat{\la}_b\tr\bSi^+\H\H^t+\log\det(\I_q+\hat{\la}_b\L^{-1})]
\end{align*}
for both the $\min(n,p)\geq r$ and the $p\geq r>n$ cases.
Hence the same arguments as in the proof of Proposition \ref{prp:SH} lead to the following proposition.
\begin{prp}\label{prp:SH_r>n}
In the model {\rm(\ref{eqn:singular_model})}, we consider the problem of estimating $\bSi$ relative to the Stein loss {\rm(\ref{eqn:singular_loss})}.
Let $c_0=6(q+1)/(3m-4q-4)$ for $3m-4q-4>0$.
Assume that $c_0q/(1+c_0)\leq b\leq C<q$ for a positive constant $C$ and $\sum_{i=1}^q\partial b/\partial \ell_i\geq 0$.
Then $\bde_r^{SH}(b)$ dominates $\bde_r^{EB}(b)$ for any possible ordering on $n$, $p$ and $r$.
\end{prp}

Further improvements on $\bde_r^{SH}(b)$ can be established in the same way as in Subsections \ref{subsec:js} and \ref{subsec:orth}.
Also, we can derive the reverse dominance relation such that $\bde_r^{EB}(b)$ dominates $\bde_r^{SH}(b)$ as in Proposition \ref{prp:SH*}.

\section{Some remarks}\label{sec:remarks}

This paper addresses the problem of estimating a high-dimensional covariance matrix of multivariate normal distribution and also discusses a unified extension to all possible cases of orderings on the dimension, the sample size and the rank of the covariance matrix.
We conclude this paper with giving some remarks.

\medskip
In this paper, it is assumed that $\bSi$ has a known rank $r$ in the singular model (\ref{eqn:singular_model}).
When $\min(n,p)\geq r$ or $p=r>n$, the observation matrix $\X$ is of rank $r$ with probability one and inherits the rank from the singular covariance matrix $\bSi$.
Thus, even if $r$ is unknown, a value of $r$ is evaluable from $\X$ as long as $\min(n,p)\geq r$ or $p=r>n$.
However the $p>r>n$ case with unknown $r$ does not permit the evaluation of $r$, which deeply affects the accuracy of estimators, particularly when $r$ is much smaller than $p$.

\medskip
Instead of the Stein loss (\ref{eqn:loss}), we may employ the quadratic loss
\begin{equation}\label{eqn:quadratic}
L_{Q}(\bde,\bSi)=\tr\bSi^{-1}(\bde-\bSi)\bSi^{-1}(\bde-\bSi).
\end{equation}
Selliah (1964) treated the $n\geq p=r$ case of covariance estimation under (\ref{eqn:quadratic}) and obtained an improved estimator based on the LU decomposition of $\S$.
For other approaches, see Haff (1979, 1980, 1991), Yang and Berger (1994) and Tsukuma (2014).
See also Konno (2009), who discussed the $p=r>n$ case under the quadratic loss (\ref{eqn:quadratic}).
For the singular case, the quadratic loss (\ref{eqn:quadratic}) probably should be replaced by
$$
L_{Q}^*(\bde,\bSi)=\tr\bSi^+(\bde-\bSi)\bSi^+(\bde-\bSi).
$$
Indeed, we can easily obtain an improved estimator similar to Selliah (1964) via the same way as in Subsection \ref{subsec:js}, but the details are omitted here.

\medskip
The observation matrix $\X$ has the form $\X=\B\Z$, where $\B$ is an unknown matrix of parameters and $\Z$ is a random matrix.
The dominance results of Section \ref{sec:unified} can be extended to the estimation problem of a scale matrix in an elliptical distribution model, where the p.d.f.\! of $\Z$ has the form $f(\tr\Z\Z^t)$ for an integrable function $f$.
The $n\geq p=r$ case with the usual Stein loss (\ref{eqn:loss}) has been studied by Kubokawa and Srivastava (1999).
Their dominance results can be extended to our singular case.

\section*{Acknowledgments}
The work is supported by Grant-in-Aid for Scientific Research (15K00055), Japan.

\end{document}